\documentclass{amsart}

\usepackage{amssymb}

\begin{document}

\title[Corrigendum IJNT 3 (2007) 119]{Corrigendum to ``On the divisibility of $\binom{n-i-1}{i-1}$ by $i$'' [Int. J. Number Theory 3 (2007) 119]}

\author{Richard J. Mathar}
\urladdr{http://www.strw.leidenuniv.nl/~mathar}
\email{mathar@strw.leidenuniv.nl}
\address{Leiden Observatory, Leiden University, P.O. Box 9513, 2300 RA Leiden, The Netherlands}


\date{\today}

\begin{abstract}
Association
of some integers $n\ge 54$
to divisibility classes
in the
Shevelev article of 2007 is corrected.
\end{abstract}

\maketitle

\section{Summary}
The examples in Section 5 of Shevelev's article on the
divisors of binomials are faulty where they
concern the cases $n=54, 60, 68, 70, 72, 78, 91$, and $96$
\cite{ShevelevIJNT3}.

\subsection{Section 1}
On page 121, the list of $n$ for which $B_n$ contains a divisor of $n$ should
be extended to 18, 45, 48, 70, 72, 75, 84, 90, 100.
We note that this list contains 74 members in the range $n\le 500$, 174 members in the range $n\le 1000$,
274 members in the range $n\le 1500$,
and 379 members in the range $n\le 2000$.

\subsection{Section 5.1}
A corrected and slightly extended list in Section 5.1 is:

$b(n)=0$. $n=$1, 2, 3, 4, 5, 6, 7, 8, 9, 10, 11, 12, 13, 15, 17, 19, 20, 21, 23, 24, 25, 29, 31, 33, 35, 37, 41, 43, 47, 49, 53, 59, 61, 67, 71, 73, 79, 83, 89, 97, 101, \ldots

$b(n)=1$. $n=$14, 16, 18, 22, 27, 28, 39, 55, 65, 77, 85, 221, 437, \ldots

$b(n)=2$. $n=$26, 30, 36, 40, 42, 44, 91, 95, 115, 119, 133, 161, 187, 247, 391, 667, 1147, \ldots

$b(n)=3$. $n=$32, 38, 45, 51, 52, 56, 57, 63, 69, 87, 145, 209, 713, 1073, \ldots

$b(n)=4$. $n=$34, 75, 84, 93, 125, 155, 203, 217, 253, 259, 299, 319, 341, 551, \ldots

$b(n)=5$. $n=$46, 48, 60, 68, 76, 81, 105, 377, 403, 407, 481, 493, 533, 589, 1189, \ldots

$b(n)=6$. $n=$50, 54, 72, 88, 99, 287, 301, 329, 527, 559, 851, \ldots

$b(n)=7$. $n=$70, 78, 80, 111, 185, 371, 451, 629, 697, 731, 799, 901, \ldots

$b(n)=8$. $n=$58, 64, 66, 74, 92, 104, 473, 517, 583, 611, 703, 779, 817, 893, 943, 1007, 1121, \ldots

$b(n)=9$. $n=$86, 117, 123, 175, 205, 343, 649, 989, 1159, \ldots

$b(n)=10$. $n=$82, 90, 96, 100, 129, 135, 140, 215, 413, 671, 689, 767, 793, 871, 1003, 1081, \ldots

$b(n)=11$. $n=$62, 120, 427, 1037, \ldots

$b(n)=12$. $n=$116, 141, 153, 235, 737, 781, 803, 869, 923, \ldots

$b(n)=13$. $n=$102, 108, 112, 136, 147, 171, 265, 949, \ldots

$b(n)=14$. $n=$106, 110, 148, 152, 295, 305, 335, 469, 913, 1139, \ldots

$b(n)=15$. $n=$94, 159, 165, 177, 183, 195, 355, 365, 1027, 1079, 1157, \ldots

$b(n)=16$. $n=$98, 124, 201, 207, 245, 395, 497, 979, 1067, 1111, \ldots

$b(n)=17$. $n=$114, 132, 164, 275, 511, 1133, \ldots

$b(n)=18$. $n=$168, 172, 189, 213, 219, 237, 249, 415, \ldots

$b(n)=19$. $n=$118, 128, 156, 184, 231, 325, 1177, \ldots

$b(n)=20$. $n=$130, 154, 261, 553, 1199, \ldots

\subsection{Section 5.2}
Associated corrections to Section 5.2 are:

$B_{39}=\{9\}$,

$B_{54}=\{4, 14, 15, 16, 20, 21\}$,

$B_{60}=\{18, 21, 22, 24, 26\}$,

$B_{68}=\{14, 22, 26, 28, 30\}$,

$B_{70}=\{4, 12, 14, 15, 22, 24, 26\}$,

$B_{72}=\{6, 14, 15, 22, 26, 33\}$,

$B_{76}=\{8, 10, 24, 26, 30\}$,

$B_{78}=\{4, 8, 10, 15, 20, 21, 28\}$,

$B_{82}=\{6, 10, 12, 14, 16, 22, 26, 28, 30, 38\}$,

$B_{91}=\{28, 35\}$,

$B_{96}=\{9, 10, 14, 20, 22, 26, 28, 38, 39, 46\}$,

$B_{98}=\{6, 10, 12, 18, 20, 22, 24, 26, 28, 30, 32, 35, 38, 42, 44, 46\}$,

Related sequences
in the OEIS are
A138389,
A178071,
A178098 -- A178101,
A178105,
A178109 and
A178110
\cite{EIS}.
\bibliographystyle{amsplain}
\bibliography{all}

\end{document}